\documentclass[12pt]{amsart}
\usepackage{amsfonts, amsmath, amsthm}

\newcommand\AAA{\mathbb{A}}
\newcommand\FF{\mathbb{F}}
\newcommand\PP{\mathbb{P}}

\newcommand\kbar{k^{\sep}}

\newtheorem{theorem}{Theorem}
\newtheorem{lemma}[theorem]{Lemma}
\newtheorem{cor}[theorem]{Corollary}

\DeclareMathOperator{\Gal}{Gal}
\DeclareMathOperator{\sep}{sep}

\begin{document}

\title[\'Etale covers of affine spaces]{\'Etale 
covers of affine spaces in positive characteristic}
\author{Kiran S. Kedlaya}
\address{Department of Mathematics, University
of California, Berkeley, Berkeley, CA 94720}
\email{kedlaya@member.ams.org}
\urladdr{math.berkeley.edu/\textasciitilde kedlaya}
\thanks{Supported by a National Science Foundation postdoctoral
fellowship.}

\date{July 17, 2002}

\begin{abstract}
We prove that every projective variety of dimension $n$
over a field of positive characteristic
admits a morphism to projective $n$-space, \'etale away from the hyperplane
$H$ at infinity, which maps a chosen divisor into $H$
and a chosen smooth point not on the divisor to some point
not in $H$.
\end{abstract}

\maketitle

\section{Introduction}

A celebrated theorem of Bely\u\i\/ \cite{belyi} asserts that a smooth,
projective, irreducible curve over the complex numbers can be defined
over a number field if and only if it admits a map to $\PP^1$ ramified
over at most three points. In positive characteristic, covers of
$\PP^1$ with even less ramification are far more prevalent: every
curve over an infinite field of characteristic $p>0$ admits a map to $\PP^1$
ramified over only one point! This assertion is both easy to prove and
surprisingly useful, especially when one wants to ``push forward'' some
problem from a complicated curve to a simple curve like the affine line.
See \cite{katz} for both the proof of the assertion (on which this note
is ultimately based)
and an application of the indicated type.

In this note, we generalize the positive characteristic assertion
to higher dimensional
varieties as follows. (Note: ``variety'' for us will mean
``separated scheme of finite type'', but not necessarily smooth, irreducible
or connected.)
\begin{theorem} \label{thm:main}
Let $X$ be a projective variety of pure dimension $n$ over an infinite
field $k$ of
characteristic $p>0$. Let $D$ be a divisor of $X$ and let $x$ be a smooth point
of $X(\kbar)$ not contained in $D$. Then there exists a morphism
$f: X \to \PP^n_k$ of $k$-schemes satisfying the following conditions:
\begin{enumerate}
\item
$f$ is \'etale away from the hyperplane $H \subseteq \PP^n$ at infinity;
\item
$f(D) \subseteq H$;
\item
$f(x) \notin H$.
\end{enumerate}
\end{theorem}
\begin{cor}
Let $X$ be a variety over $k$, and let $x$ be a smooth point
of $X(\kbar)$. Then there exists a finite \'etale morphism $f: U \to \AAA^n$
for $U$ some open dense subset of $X$, defined over $k$ and containing $x$.
\end{cor}
By noetherian induction, $X$ can thus be covered with open (necessarily
affine) subsets which are finite \'etale covers of $\AAA^n$.

Note that the theorem is strictly stronger than the corollary; from the
corollary, one only gets a rational map from $X$ to $\PP^n$. However, in
some cases it may be the corollary that is most directly useful, again
when one needs to ``push forward'' a problem to a simpler space via an
\'etale map, but only on an open dense subset of the original space.
One example of this situation is the author's proof of finite dimensionality of
rigid cohomology with coefficients \cite{me:finite}.

The restriction to infinite $k$ is probably not necessary; it intervenes in
the proof because we must choose certain constructions ``generically'' to
avoid undesired behaviors. In any case, if $k$ is finite, the conclusion of
the theorem holds over some finite extension of $k$, depending on the rest
of the input data.

\section{Proof of the Theorem}

To prove the theorem, we will need to string together a chain of carefully
chosen maps. To facilitate this, we make the following definitions.
A \emph{good triple} will always mean a triple $(Y, E, y)$, where
$Y$ is a variety of pure dimension $n$ over $k$, $E$ is a divisor of $Y$
defined over $k$, and $y$ is a point of $Y(\kbar)$ not contained in $E$.
Given two good triples $(Y_1, E_1, y_1)$ and $(Y_2, E_2, y_2)$, a 
\emph{good morphism} will be a finite morphism $f: Y_1 \to Y_2$ of $k$-schemes,
with $f(E_1) \subseteq E_2$, $f(y_1) = y_2$, and $f$ \'etale on
$Y_2 \setminus E_2$.

In this language, the given triple $(X, D, x)$ is good, and the
problem is to find a chain of good morphisms leading from $(X,D,x)$
to $(\PP^n, H, z)$ for some $z \notin H$. We construct this chain in three
steps.

Reminder: the assertion ``property X holds for the generic Y'' means that
property X holds for all Y's in an open dense subset of the natural parameter
space of all objects Y. In particular, this type of assertion is stable
under conjuction on property X.

\subsection*{Step 1: Noether normalization}

For our first step, we construct a good morphism $\pi: (X, D, x) \to
(\PP^n, D_0, x_0)$ by Noether normalization. Choose a projective
embedding $g: X \to \PP^m$ of $X$.
For a generic $(m-n-1)$-plane $P$ in $\PP^m$, the map $\pi: X \to \PP^n$
induced by projection away from $P$ is finite and
has the following additional properties:
\begin{enumerate}
\item[(a)]
$\pi(x) \notin \pi(D)$, that is, $P$ does not meet the join $J$ of $x$
and $D$. That is because $\dim J + \dim P = n + (m-n-1) < m$.
\item[(b)]
$\pi$ is \'etale over $\pi(x)$. This follows from Bertini's theorem
and the fact that a generic $(m-n)$-plane through $x$
is the intersection of $n-m$ generic hyperplanes:
the intersection of $X$ with one generic hyperplane is smooth, the
intersection of the result with a second generic hyperplane is again
smooth, and so on, until the intersection of $X$ with the $(m-n)$-plane
is smooth and hence reduced.
\end{enumerate}
Fixing a choice of $P$, take
$x_0 = \pi(x)$ and $D_0$ to be the union of $\pi(D)$ with the branch 
locus of $\pi$.

We have now eliminated all of the geometry of the ambient variety $X$
from the discussion; the rest of the argument takes place within
the projective space $\PP^n$.

\subsection*{Step 2: Additive polynomials}

In this step, we construct a sequence of good triples $(\PP^n,
x_i, D_i)$ for $i=0, \dots, n$, starting with the good triple
$(\PP^n, x_0, D_0)$ from the previous step, and a sequence of
good morphisms
$f_i: (\PP^n, x_{i}, D_{i}) \to (\PP^n, x_{i+1}, D_{i+1})$,
such that $D_i$ is the union of $i$ hyperplanes $H_{ij}$ ($j=0,
\dots, i-1$) in general position (that is, whose mutual intersection
has codimension $i$) with the cone $C_i$ over a hypersurface within
a plane of codimension
$n-i-1$ in $\PP^n$ not meeting $\cap_j H_{ij}$. In particular, $D_n$
will be the union of $n+1$ hyperplanes meeting transversely.

Before proceeding to the construction, we recall a bit of algebra
peculiar to positive characteristic.
\begin{lemma}
Let $R$ be a ring of characteristic $p>0$. Then for any polynomial
$P \in R[t]$ of degree $m$, there is a canonical multiple $Q$ of $P$ having
the form
\[
Q(t) = \sum_{i=0}^m r_i t^{p^i}
\]
for some $r_i \in R$ with $r_m$ nonzero.
Moreover, if $R = k[x_1, \dots, x_l]$ and $P$ is homogeneous as a polynomial
in $x_1, \dots, x_l, t$, then so is $Q$.
\end{lemma}
A polynomial of the form prescribed for $Q$ is called \emph{additive}, since
such polynomials are precisely those for which $Q(t+u) = Q(t)+Q(u)$
identically. The proof of the lemma is standard, but as it may not be familiar
to all readers, we include it.
\begin{proof}
In the ring $\FF_p[t_1, \dots, t_m, x]$, define $s_1, \dots, s_m$ 
as the elementary symmetric functions of $t_1, \dots, t_m$:
\[
x^m + s_1 x^{m-1} + \cdots + s_m = (x+t_1)\cdots (x+t_m).
\]
Then the polynomial
\[
S(x) =
\prod_{h_1, \dots, h_m \in \FF_p} (x + h_1 t_1 + \cdots + h_m t_m)
\]
is symmetric in $t_1, \dots, t_m$,
so its coefficients can be expressed
as polynomials in the $s_i$.
On the other hand, up to sign, $S(x)$ is the value of the Moore determinant
\[
\det \begin{pmatrix}
x & x^p & \cdots& x^{p^m} \\
t_1 & t_1^p & & t_1^{p^m} \\
\vdots & \ddots & \\
t_m & t_m^p &  & t_m^{p^m}
\end{pmatrix},
\]
by the same argument used to evaluate the Vandermonde determinant: the
determinant clearly vanishes when any one of the linear factors of $S$ is
set to zero, and has the same total degree as $S$. Thus $S$ is additive
in $x$.
The desired $Q$ is now the image of $S$ under the homomorphism
from $\FF_p[s_1, \dots, s_m, x]$ to $R[t]$ sending $x$ to $t$ and
$S_i$ to the coefficient of $t^{m-i}$ in $P$.
\end{proof}

We first outline the construction of $f_i$ given $x_{i}$ and
$D_{i}$, then record the geometric conditions that
must be satisfied for the construction to go through. The construction
will depend on a choice of homogeneous coordinates $z_0, \dots, z_n$
such that $H_{ij}$ is the zero locus of $z_j$ and the defining equation $P_i$
of $C_i$ depends only on $z_i, \dots, z_n$. The set of such choices forms
an irreducible parameter variety; we will ultimately show that each
of the necessary conditions is satisfied on an open and nonempty, so
dense, subset of the parameter variety.

Regard $P_i$ as a polynomial in $z_i$ whose coefficients are polynomials 
in $z_{i+1}, \dots, z_n$. Let $Q_i$ be the multiple of $P_i$ produced by
the previous lemma, and put $d_i = \deg(Q_i)$.
Now define the map $f_i:
\PP^n \to \PP^n$ by sending $(z_0:\cdots:z_n)$ to $(w_0:\cdots:w_n)$, where
\[
w_j = \begin{cases}
z_j^{d_i} - z_j z_n^{d_i-1} & j \neq i, n \\
Q_i(z_i, \dots, z_n) & j = i \\
z_n^{d_i} & j = n
\end{cases}
\]
and take $x_{i+1} = f_i(x_{i})$, $H_{(i+1)j}$ to be the zero locus of $w_j$
for $j=0, \dots, i$, and $C_{i+1}$ to be the zero locus of the constant
coefficient of $Q_i$. In particular, $C_{i+1}$ is the zero locus of a 
polynomial depending only on $z_{i+1}, \dots, z_n$.

For $f_i$ to be a regular map, the $w_j$ must have no common zeroes.
In that case, the non\'etale locus of $f_i$ is contained in the zero 
locus of $z_n$ times the constant coefficient of $Q_i$. 
In short, the construction gives what we want provided that
the following conditions hold.
\begin{enumerate}
\item[(a)]
The degree of $P_i$ as a polynomial in $z_i$ alone is equal to its total
degree in $z_i, \dots, z_n$.
\item[(b)]
$z_n$ and $z_j^{d_i} - z_j z_n^{d_i-1}$ take nonzero values
at $x_i$ for $j \neq i, n$.
\item[(c)]
$Q_i$ has nonzero constant coefficient.
\item[(d)]
$Q_i$ takes a nonzero value at $x_i$.
\end{enumerate}
Each of these is clearly an open condition on the parameter variety of
coordinate systems $z_0, \dots, z_n$. We conclude the construction by
verifying that each condition is not identically violated. Then each 
condition holds on an open dense subset of the parameter variety;
since $k$ is infinite, the intersection of these open dense subsets
contains infinitely many $k$-rational points, any one of which
yields a satisfactory choice of $f_i$.

The first two conditions are clearly not identically violated.
To check (c) and (d), it suffices to work in the
projection from $\cap H_{ij}$. In the image of this projection,
draw the line through the image of $x_i$ and the point with $z_0=\cdots=z_{n-1} 
= 0$, and choose an identification of this line $\PP^1$ in which
the latter point becomes $\infty$. Then (c) is satisfied
if the intersections of $C_i$ with this line are identified with a set
of elements of $k^{\sep}$ which are linearly independent over $\FF_p$
(and in that case $d_i = p^{\deg P_i}$),
and (d) is satisfied if the same holds after including
$x_i$ as well.

We now turn the tables, regarding the line as fixed and varying the coordinate
system, under the constraint that the point with $z_0=\cdots=z_{n-1}=0$
remains on the line. As we do this, the elements of $k^{\sep}$ that we wrote
down previously are moved around by linear fractional transformations,
and by the following lemma, at some point they become linearly
independent over $\FF_p$.

\begin{lemma}
Let $\{r_1, \dots, r_m\}$ be a finite subset of $\kbar$ stable under
$\Gal(\kbar/k)$. Then for a generic choice of $a,b,c,d \in k$
(i.e., away from a Zariski closed subset of $\AAA^4_k$), 
if we set $\tau(x) = (a+bx)/(c+dx)$, then
\[ h_1 \tau(r_1) + \cdots + h_m \tau(r_m) \neq 0
\]
for any $h_1, \dots, h_m \in \FF_p$ not all zero.
\end{lemma}
\begin{proof}
It suffices to check this separately for each choice of
$h_1, \dots, h_m$, since there are finitely many such choices.
Moreover, it is enough to check this under the additional restriction
$a=b=0$ and $d=1$. In that case, the expression in question becomes
\[
\frac{h_1}{c + x_1} + \cdots + \frac{h_m}{c+x_m} = \frac{R'(-c)}{R(c)},
\]
where $R(x) = \prod_j (x + x_j)^{h_j}$. Since $R(x)$ is not a $p$-th
power, its derivative does not vanish identically. Thus the expression
does not vanish identically over all choices of $a,b,c,d$,
as desired.
\end{proof}

Thus (c) and (d) hold for some coordinate system, completing the
verification of the necessary conditions for the construction of $f_i$.

\subsection*{Step 3: The Abhyankar map}

For the third step, we must construct a good morphism $f_n: (\PP^n, x_n, D_n)
\to (\PP^n, x_{n+1}, H)$ for some $x_{n+1}$, where $D_n$ is the union of
$n$ transverse hyperplanes. We explicitly construct this morphism
by writing down polynomials $g_i$
in the variables $z_0,\dots, z_n$, for $i=0, \dots, n$, as follows.
For each $(i+1)$-element subset $I = \{j_0, \dots, j_i\}$
of $\{0, \dots, n\}$, with $j_0 < \cdots < j_i$, define
\[
m_I = z_{j_0}^{1+p+\cdots+p^{n-i}}
z_{j_1}^{p^{n-i+1}} \cdots z_{j_i}^{p^n}
\]
and let $g_i$ be the sum of the $m_I$ over all $(i+1)$-element subsets $I$.
For example, when $n=2$, we have
\begin{align*}
g_0 &= z_0^{p^2+p+1} + z_1^{p^2+p+1} + z_2^{p^2+p+1} \\
g_1 &= z_0^{p+1} z_1^{p^2} + z_0^{p+1} z_2^{p^2} + z_1^{p+1} z_2^{p^2} \\
g_2 &= z_0 z_1^p z_2^{p^2}.
\end{align*}

Let us observe some facts about the $g_i$. First, they are all homogeneous
of degree $1+p+\cdots+p^n$. Second, they have no common zero
except $z_0 = \cdots = z_n = 0$, by the same argument as for the
elementary symmetric functions: if $g_n = 0$, then one of the $z_i$ must
be zero; in that case, if $g_{n-1} = 0$, then another of the $z_i$
must be zero, and so on.
These two facts allow us to define a morphism $f_n: \PP^n \to \PP^n$ by
the formula
\[
(z_0:\cdots:z_n) \mapsto
(g_0:\cdots:g_n).
\]
Third, note that for $z_0, \dots, z_n$ all nonzero, the differentials
$dg_0, \dots, dg_n$ are linearly independent. Namely, $dg_n$ is a nonzero
multiple of $dz_0$; $dg_{n-1}$ is a nonzero multiple of $dz_1$ plus a 
multiple of $dz_0$; $dg_{n-2}$ is a nonzero multiple of $dz_2$ plus a
linear combination of $dz_0$ and $dz_1$; and so on. This means that $f_n$
is \'etale away from the zero locus of $z_0\cdots z_n$, i.e.,
the zero locus of $g_n$. 

In passing, we note that
the case $n=1$ of this construction yields what is commonly called the
Abhyankar map, which expresses the affine line minus a point as an \'etale
cover of the full affine line. It seems a fitting tribute to Abhyankar's
work to bestow the same name on this higher-dimensional analogue.

To conclude, if we set
$x_{n+1} = f_n(x_n)$, and $H$ equal to the hyperplane
$z_{n} = 0$, the map $f_n$ gives a good morphism from
$(\PP^n, x_n, D_n)$ to $(\PP^n, x_{n+1}, H)$. Stringing together
the good morphisms $f_0, \dots, f_n$ yields a good morphism from
$(X, x, D)$ to $(\PP^n, x_{n+1}, H)$, completing the proof of the theorem.

\section{Remarks}

To conclude, we speculate briefly about higher-dimensional generalizations
of the original theorem of Bely\u\i. One such is a theorem of Bogomolov
and Pantev \cite{bp}: given a variety $X$ over an algebraically closed field
of characteristic zero and a proper subset $D$ of $X$, there exists a 
toroidal embedding $E \subset Y$ (in the sense of \cite{kkms}), a blowup
$\tilde{X}$ of $X$, and a morphism
$f: X \to Y$, \'etale away from $E$, with $f(\tilde{D}) \subseteq E$.
This theorem
leads instantly to a proof of weak resolution of singularities in positive
characteristic, as $\tilde{E} \subset \tilde{X}$ is then also a toroidal
embedding, to which the toroidal resolution algorithm of \cite{kkms}
may be applied. (It must be noted that the final conclusion, that $\tilde{E}
\subset \tilde{X}$ is a toroidal embedding because it is a finite \'etale
cover of another toroidal embedding,
fails severely in positive characteristic, so it is unclear
whether the Bogomolov-Pantev approach can be used to say anything about
resolution of singularities there.)

The Bogomolov-Pantev result is quite beautiful, but one can ask whether it
is possible to formulate a higher-dimensional statement more in the spirit
of the theorem of Bely\u\i. Such a statement would, for one, apply to
a restricted class of $X$ but would give a conclusion about a map from $X$
itself and not some blowup. It would also have a specific target space,
say a projective space missing some hyperplanes; the Bogomolov-Pantev target
space $Y$
ends up being a rational variety, but besides that one can say little about
it.

One natural generalization of three points in projective space is $n+2$
hyperplanes in projective $n$-space. With this, we formulate the following
question: for which varieties $X$ and divisors $D$ over a number field
does there exist
a morphism $f: X \to \PP^n$ such that $f(D)$ is contained in the union $S$
of $n+2$ hyperplanes of $\PP^n$, and $f$ is \'etale over $\PP^n \setminus S$?
Again, the pair $D \subset X$ must be a toroidal embedding, but are there
other conditions?

\end{document}